\input amstex
\input amsppt.sty
\magnification=\magstep1
\baselineskip=16truept
\vsize=22.2truecm
\nologo
\pageno=1
\topmatter

\def\pmod #1{\ (\roman{mod}\ #1)}
\def\Proof{\noindent{\it Proof}}

\title
A $q$-analogue of Wolstenholme's harmonic series congruence
\endtitle
\author
Ling-Ling Shi$^1$ and Hao Pan$^2$
\endauthor
\address
1. College of Mathematics and Physics, Zhejiang Normal University,
Jinhua 321004, People's Republic of China
\endaddress
\email{linglingshi\@eyou.com}\endemail
\address
2. Department of Mathematics, Nanjing University,
Nanjing 210093, People's Republic of China
\endaddress
\email{haopan79\@yahoo.com.cn}\endemail
\endtopmatter
\TagsOnRight
\document
The $n$-th harmonic number $H_n$ is defined by
$$
H_n=\sum_{j=1}^n\frac{1}{j}.
$$
A well-known result of Wolstenholme [3] asserts that if $p\geq 5$ is prime, then
$$
H_{p-1}\equiv 0\pmod{p^2}.\tag 1
$$
We define the $q$-harmonic numbers by
$$
H_n(q)=\sum_{j=1}^n\frac{1}{[j]_q},
$$
where
$$
[n]_q=(1-q^n)/(1-q)=1+q+\cdots+q^{n-1}.
$$
In view of the $q$-analogue of Glaisher's congruence, Andrews ([1], Theorem 4) showed that
$$
H_{p-1}(q)\equiv \frac{p-1}{2}(1-q)\pmod{[p]_q}\tag 2
$$
for any odd prime $p$. In the present paper, we will extend Andrews' result as follows.
\proclaim{Theorem 1} Let $p\geq 5$ be a prime. Then
$$
H_{p-1}(q)\equiv \frac{p-1}{2}(1-q)+\frac{p^2-1}{24}(1-q)^2[p]_q\pmod{[p]_q^2}.\tag 3
$$
\endproclaim
The following lemma is a $q$-analogue of the congruence
$$
\sum_{j=1}^{p-1}\frac{1}{j^2}\equiv 0\pmod{p}
$$
for prime $p\geq 5$.
\proclaim{Lemma 2} For any prime $p\geq 5$,
$$
\sum_{j=1}^{p-1}\frac{q^j}{[j]_q^2}\equiv-\frac{p^2-1}{12}(1-q)^2\pmod{[p]_q}\tag 4
$$
and
$$
\sum_{j=1}^{p-1}\frac{1}{[j]_q^2}\equiv-\frac{(p-1)(p-5)}{12}(1-q)^2\pmod{[p]_q}.\tag 5
$$
\endproclaim
\Proof. Clearly
$$
\sum_{j=1}^{p-1}\frac{q^j}{[j]_q^2}=(1-q)^2\sum_{j=1}^{p-1}\frac{q^j}{(1-q^j)^2}.
$$
We set
$$
G(q)=\sum_{j=1}^{p-1}\frac{q^j}{(1-q^j)^2}.
$$
Let $\zeta=e^{2\pi i/p}$. Since
$$
[p]_q=\frac{1-q^p}{1-q}=\prod_{m=1}^{p-1}(q-\zeta^m),
$$
in order to prove $(4)$, we only need to show that
$$
G(\zeta^m)=-\frac{p^2-1}{12}
$$
for $m=1,2,\ldots,p-1$. As $p$ is prime, we have
$$
G(\zeta^m)=\sum_{j=1}^{p-1}\frac{\zeta^{mj}}{(1-\zeta^{mj})^2}=\sum_{j=1}^{p-1}\frac{\zeta^j}{(1-\zeta^j)^2}=G(\zeta)
$$
provided that $1\leq m<p$. Define
$$
G(q,z)=\sum_{j=1}^{p-1}\frac{q^j}{(1-q^jz)^2}.
$$
Observe that
$$
\sum_{j=1}^{p-1}\zeta^{kj}=\cases p-1\qquad&\text{if }p\mid k,\\-1\qquad&\text{if }p\nmid k.\endcases
$$
Thus when $|z|<1$,
$$
\aligned
G(\zeta,z)=&\sum_{j=1}^{p-1}\zeta^j\sum_{k=0}^{\infty}\zeta^{kj}(k+1)z^k\\
=&\sum_{k=1}^{\infty}kz^{k-1}\sum_{j=1}^{p-1}\zeta^{kj}\\
=&p\sum_{k=1}^{\infty}(pk)z^{pk-1}-\sum_{k=1}^{\infty}kz^{k-1}\\
=&\frac{p^2z^{p-1}}{(1-z^p)^2}-\frac{1}{(1-z)^2}.
\endaligned
$$
Letting $z\to 1$ and using L'Hospital's rule, we obtain that
$$
G(\zeta)=G(\zeta,1)=\lim_{z\to 1}\frac{p^2z^{p-1}(1-z)^2-(1-z^p)^2}{(1-z^p)^2(1-z)^2}=\frac{1-p^2}{12}
$$
as desired.

Finally, since
$$
\sum_{j=1}^{p-1}\frac{1}{[j]_q^2}=\sum_{j=1}^{p-1}\frac{1-q^j+q^j}{[j]_q^2}=(1-q)\sum_{j=1}^{p-1}\frac{1}{[j]_q}+\sum_{j=1}^{p-1}\frac{q^j}{[j]_q^2},
$$
(5) immediately follows from (2) and (4).
\qed
{\medskip\noindent{\it  Proof of Theorem 1}}.
$$
\aligned
&H_{p-1}(q)-\frac{p-1}{2}(1-q)\\
=&\frac{1-q}{2}\sum_{k=1}^{p-1}\big(\frac{1}{1-q^k}+\frac{1}{1-q^{p-k}}-1\big)\\
=&\frac{1-q}{2}\sum_{k=1}^{p-1}\frac{1-q^p}{(1-q^k)(1-q^{p-k})}.
\endaligned
$$
Then by (4) we have
$$
\aligned
&\frac{(1-q)(1-q^p)}{2}\sum_{k=1}^{p-1}\frac{1}{(1-q^k)(1-q^{p-k})}\\
=&\frac{(1-q)(1-q^p)}{2}\sum_{k=1}^{p-1}\frac{q^k}{(1-q^k)(q^k-q^p)}\\
\equiv&-\frac{(1-q)(1-q^p)}{2}\sum_{k=1}^{p-1}\frac{q^k}{(1-q^k)^2}\\
\equiv&\frac{p^2-1}{24}(1-q)(1-q^p)\pmod{[p]_q^2}.
\endaligned
$$
This concludes the proof of $(3)$.\qed

\enddocument